\documentclass[a4paper,12pt]{article}
\usepackage[top=2.0cm,bottom=3.0cm,left=2.5cm,right=2.5cm]{geometry}
\usepackage{fancyhdr,graphicx,amsfonts,amssymb,amsthm,newlfont,amsmath}
\usepackage{bm,latexsym,mathrsfs,setspace,lipsum}

\begin{document}

\title{{The size of $k$-th order generalized Fibonacci cubes
\thanks{This work is supported by Natural Science Foundation of Shandong Province(ZR2025MS62),
National Natural Science Foundation of China (No. 12171414)
and Taishan Scholars Special Project of Shandong Province.}}}

\author{Jianxin Wei$^{a}$ \footnote{Corresponding author.},
Yujun Yang$^{b}$\\
\scriptsize{$^{a}$School of Mathematics and Statistics Science, Ludong University, Yantai, Shandong, 264025, PR China}\\
\scriptsize{$^{b}$School of Mathematics and Information Science, Yantai University, Yantai, Shandong, 264005, PR China}\\
\scriptsize{E-mail address:
wjx0426@ldu.edu.cn,
yangyj@ytu.edu.cn}}
\date{}
\maketitle
\textbf{Abstract:}
Let $k\geq2$.
Then the $k$-th order Fibonacci cube $\Gamma^{(k)}_{n}$ is the subgraph of the hypercube $Q_{n}$
induced by vertices without $k$ consecutive $1$s.
The case $k=2$ corresponds to the classic Fibonacci cube $\Gamma_{n}$.
There are three kinds of calculation formulas of the size of $\Gamma_{n}$:
the iteration form $|E(\Gamma_{n})|=|E(\Gamma_{n-1})|+|E(\Gamma_{n-2})|+F_{n}$ (Hsu, 1993), 
the convolution form $|E(\Gamma_{n})|=\mathop{\sum}\limits_{i=1}^{n}F_{i}F_{n-i+1}$ (Klav\v{z}ar, 2005) 
and the linear form $|E(\Gamma_{n})|=\frac{nF_{n+1}+2(n+1)F_{n}}{5}$ (Munarini et al., 2001). 
Belbachir and Ould-Mohamed (2020) studied the iteration and convolution formulas of the size of $\Gamma^{(3)}_{n}$.
Very recently,
Mollard (2025) deduced the iteration formula of the size of $\Gamma^{(k)}_{n}$ for $k\geq2$.
In this paper,
we give the the formulas of convolution and linear forms of $|E(\Gamma^{(k)}_{n})|$ for all $k\geq2$.
Specifically,
we obtain the formula of $|E(\Gamma^{(k)}_{n})|$ in terms of convolved $k$-th order Fibonacci numbers and
the formula of $|E(\Gamma^{(k)}_{n})|$ of linear expression of $k$ consecutive $k$-th order Fibonacci numbers.

\newcommand{\trou}{\vspace{1.5 mm}}
\newcommand{\noi}{\noindent}
\newcommand{\ol}{\overline}

{\vspace{1.5 mm}}

\section{Introduction}

The \emph{hypercube of dimension} $n$,
denoted by $Q_{n}$,
is the $n$-regular graph of $2^{n}$ vertices labelled by binary strings of length $n$,
an edge joining two vertices whenever the corresponding strings differ in exactly one coordinate.
It is known that the size of $Q_{n}$ is $n2^{n-1}$.

The \emph{Fibonacci cube} $\Gamma_{n}$ was introduced by Hsu \cite{Hsu} as the subgraph of $Q_{n}$ induced by the binary strings that contain no two consecutive 1s.
Fibonacci cube has many attractive structural,
numerous combinatorial and beautiful properties,
and many applications.
For more information on Fibonacci cubes,
see the survey \cite{Klavzar} and the recent book \cite{EgeciogluKlavzarMollard}.

Inspired by the fact that
the vertex set of $\Gamma_{n}$ is consisted of all binary strings of length $n$ in which `$11$' is forbidden,
for a given binary string $f$,
Ili\'{c},
Klav\v{z}ar
and Rho \cite{IlicKlavzarRho1} defined
the \emph{generalized Fibonacci cube} $Q_{n}(f)$ as the graph obtained from
the hypercube $Q_{n}$ by removing all vertices that contain $f$ as a factor.
Note that the term ``generalized Fibonacci cubes'' has been used in \cite{HsuChung}
for a restricted family of the graph $Q_{n}(1^{k})$ for $k\geq2$.
Following Mollard \cite{Mollard} the graph $Q_{n}(1^{k})$ is called \emph{$k$-th order Fibonacci cube}
and the notation $\Gamma_{n}^{(k)}$ is used for it in this paper.
It is obvious that the case $k=2$ corresponds to the Fibonacci cube $\Gamma_{n}$.
The graphs $\Gamma_{n}^{(3)}$ ($n\in[5]$) are shown in Fig. 1.
\begin{center}
\includegraphics[scale=0.75]{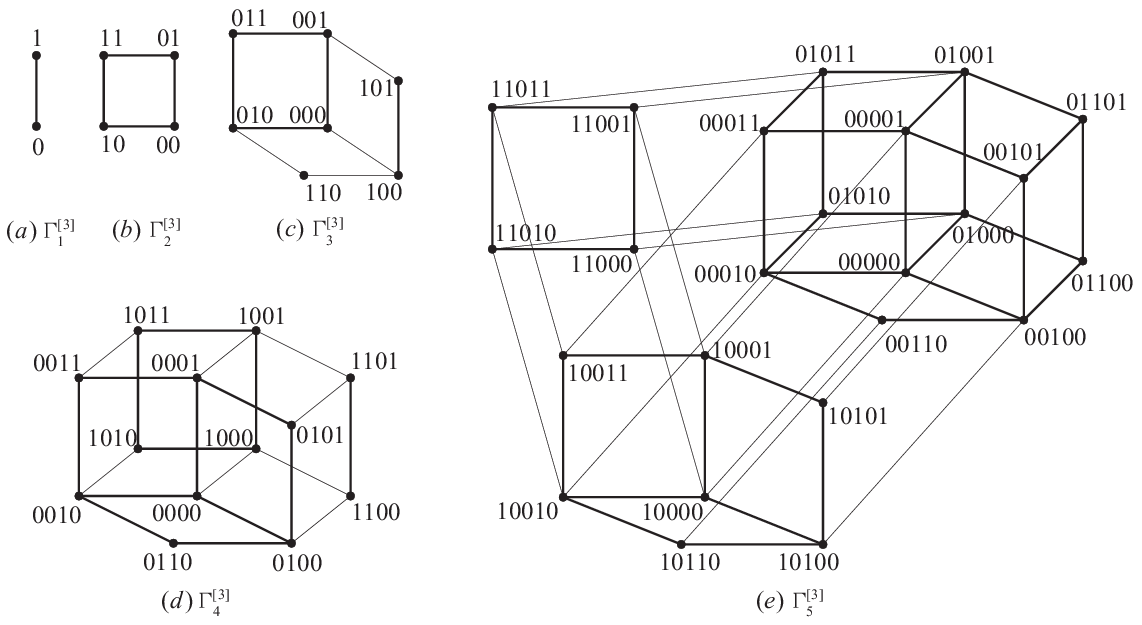}\\
{\footnotesize Fig. 1. The graphs $\Gamma_{n}^{(3)}$, $n\in[5]$.}
\end{center}

There have been some studies on the $k$-th order Fibonacci cubes.
In \cite{HsuChung},
the diameter of $\Gamma_{n}^{(k)}$ was determined,
and it was shown that the distance between the vertices of $\Gamma_{n}^{(k)}$ is their Hamming distance.
In \cite{LiuHsuChung},
the Hamiltonicity of $\Gamma_{n}^{(k)}$ was characterized.
In \cite{Salvi},
it was proven that every edge of $\Gamma_{n}^{(k)}$,
but few initial cases,
belongs to cycles of every even length.
In \cite{KlavzarRho},
the Wiener index of $\Gamma_{n}^{(k)}$ was studied.
The $3$-th order Fibonacci cube $\Gamma_{n}^{(3)}$ cube was called
Tribonacci cube by Belbachir and Ould-Mohamed and some enumerative
properties and cube polynomial of $\Gamma_{n}^{[3]}$ were studied \cite{BelbachirOuldMohamed}.
In \cite{RhoVesel},
a linear time algorithm for the recognition of Tribonacci cubes was given.
In \cite{WangNiu},
the Mostar index of Tribonacci cubes were calculated.
In \cite{Mollard},
for general $k\geq2$,
the order, size and cube polynomial of $\Gamma_{n}^{(k)}$ as well as their generating functions were studied.

We pay a special attention on the size of the $k$-th order Fibonacci cube.
We find that there are three kinds of calculation formulas of the size of $\Gamma_{n}$.
In the seminal paper \cite{Hsu},
Hus gave the iterative formula of the size of $\Gamma_{n}$ depend on the decomposition of $\Gamma_{n}$:

$$|E(\Gamma_{n})|=|E(\Gamma_{n-1})|+|E(\Gamma_{n-2})|+F_{n},n\geq2.$$
\\
Klav\v{z}ar \cite{Klavzar1} expressed the formula of the size of $\Gamma_{n}$ as the convolution product of Fibonacci numbers:
$$|E(\Gamma_{n})|=\mathop{\sum}\limits_{i=1}^{n}F_{i}F_{n-i+1},n\geq1.$$
\\
Munarini et al. \cite{MunariniCippoSalvi} provided the formula of
$|E(\Gamma_{n})|$ expressed linearly by two consecutive Fibonacci numbers:
$$|E(\Gamma_{n})|=\frac{nF_{n+1}+2(n+1)F_{n}}{5},n\geq1.$$
\\
Belbachir and Ould-Mohamed \cite{BelbachirOuldMohamed} studied the size of $\Gamma^{(3)}_{n}$.
First they obtained the recurrence relation of $|E(\Gamma^{(3)}_{n})|$ by the decomposition of $\Gamma^{(3)}_{n}$,
and then they shown that the size of $\Gamma^{(3)}_{n}$ is the convolution product of the $k$-th order Fibonacci numbers depending on the generating function of $|E(\Gamma^{(3)}_{n})|$.
\\
Mollard \cite{Mollard} recently deduced the iteration formula of $|E(\Gamma^{(k)}_{n})|$ and the generating function of $|E(\Gamma^{(k)}_{n})|$ for $k\geq2$.

Those studies on size of $\Gamma^{(k)}_{n}$ naturally lead to the following question:
are there the convolution and linear formulas of $|E(\Gamma^{(k)}_{n})|$ for general $k$ $(\geq2)$?
In this paper,
we give a positive answer to this question.
More specifically,
we show that $|E(\Gamma^{(k)}_{n})|$ is the convolution product of $k$-th order Fibonacci numbers,
and $|E(\Gamma^{(k)}_{n})|$ also is the linear combination of $k$ consecutive $k$-th order Fibonacci numbers.

This paper is organized as follows.
In Section 2,
the $k$-th order Fibonacci numbers and some results about $\Gamma^{(k)}_{n}$ are recalled.
In Section 3,
the formulas of convolution and linear forms of $|E(\Gamma^{(k)}_{n})|$ are determined.
In the last section,
we suggest to give the convolution and linear formulas for more Fibonacci-like cubes,
and as an example,
the convolution formula of Fibonacci $p$-cubes is determined.

\section{Preliminaries}

In this section,
we give some concepts,
notations and results needed this paper.
To avoid ambiguity with initial conditions,
we first introduce the Fibonacci numbers.
The well known \emph{Fibonacci numbers} are given by the initial values $F_{0}=0$ and $F_{1}=1$ and recurrence
$$F_{n}=F_{n-1}+F_{n-2}, n\geq2.$$

Fibonacci numbers spawn many generalizations.
A classic one of them is to consider the sequence in which each element is the sum of the previous $k$ elements.
For a given integer $k\geq2$,
the \emph{$k$-th order Fibonacci numbers} \cite{Flores,Miles} (also known as the \emph{Fibonacci $k$-step numbers})
can be defined by the following the initial values and recurrence
$$F_{0}^{(k)}=\ldots=F_{k-2}^{(k)}=0, F_{k-1}^{(k)}=1,
F_{n}^{(k)}=F_{n-1}^{(k)}+F_{n-2}^{(k)}+\ldots+F_{n-k}^{(k)}, n\geq k.~~~~~~~~~~~~~~~~~~~~~~~~~~~~(2.1)$$

Obviously,
$F_{n}^{(2)}$ is the usual Fibonacci number $F_{n}$.
The $k$-th generalized Fibonacci numbers are shown in Table 1 for $k=2,3,4$ and $5$.
\begin{table}[!htbp]
\footnotesize
\caption{$F^{(k)}_{n}$, $k=2,3,4,5$ and $n=0,1,2,\ldots,16$.}
\centering
  \begin{tabular}{c|c|c|c|c|c|c|c|c|c|c|c|c|c|c|c|c|c}
  \hline
  $n$             &0& 1& 2& 3& 4& 5& 6& 7& 8& 9& 10& 11& 12& 13& 14& 15& 16 \\
  \hline
   $F^{(2)}_{n}$ & 0& 1&1 & 2&3 & 5&8 & 13& 21&34&55 &89 & 144& 233&377 &610 &987 \\
  \hline
  $F^{(3)}_{n}$ & 0&0 &1 &1 &2 &4 &7 &13 & 24& 44& 81& 149& 274& 504&927&1705&3136\\
  \hline
  $F^{(4)}_{n}$ & 0&0 &0&1 &1 &2 &4 & 8&15& 29& 56& 108& 208&401& 773& 1490&2872 \\
  \hline
  $F^{(5)}_{n}$ & 0&0 &0&0 &1 &1 &2 & 4&8& 16& 31& 61& 120& 236& 464&912&1793 \\
  \hline
  \end{tabular}
\end{table}

Recall that $\Gamma_{n}^{(k)}$ is the subgraph of $Q_{n}$ induced by the vertices in which the substring $1^{k}$ is forbidden.
For convenience,
the vertex set $V(\Gamma_{n}^{(k)})$ of $\Gamma_{n}^{(k)}$ is denoted by $\mathbf{F}_{n}^{(k)}$.
In other words,
$\mathbf{F}_{n}^{(k)}$ is the set of strings of length $n$ with at most $k-1$ consecutive 1s.
It can be seen that if $n\leq k-1$ for a given $k$,
then  $\Gamma_{n}^{(k)}\cong Q_{n}$ and so $\mathbf{F}_{n}^{(k)}$ is the set of all binary strings of length $n$.
For $n\geq k$,
the researches on the $k$-th order Fibonacci cubes $\Gamma_{n}^{(k)}$ are based on the following partition of $\mathbf{F}_{n}^{(k)}$.

\trou \noi {\bf Proposition 2.1} \cite{HsuChung,Mollard}.
\emph{Let $k\geq2$ and $n\geq k$. Then}
$\mathbf{F}_{n}^{(k)}=\mathop{\bigcup}\limits_{i=0}^{k-1}1^{i}0\mathbf{F}_{n-i-1}^{(k)}.$

It is well-known that the order of $\Gamma_{n}$ is $F_{n+2}$.
From the above analysis and Proposition 2.1,
the following result holds on the $k$-th order Fibonacci cubes.

\trou \noi {\bf Proposition 2.2} \cite{HsuChung,Mollard}.
\emph{Let $k\geq2$ and $n\geq1$.
Then the order of $\Gamma_{n}^{(k)}$ is $F_{n+k}^{(k)}$.}

By considering the edges between $1^{i}0\mathbf{F}_{n}^{(k)}$
and $1^{j}0\mathbf{F}_{n}^{(k)}$ for $0\leq i\neq j\leq k-1$,
the fundamental decomposition of $\Gamma_{n}^{(k)}$ can be obtained as shown in Fig. 2.

Based on the decomposition,
Mollard deduced the following recurrence relation on the size of $k$-th order Fibonacci cubes.

\trou \noi {\bf Theorem 2.3} \cite{Mollard}.
\emph{Let $n\geq k\geq2$.
Then}
$|E(\Gamma_{n}^{(k)})|=\mathop{\sum}\limits_{i=1}^{k}|E(\Gamma_{n-i}^{(k)})|
+\mathop{\sum}\limits_{i=2}^{k}(i-1)|V(\Gamma_{n-i}^{(k)})|.$

For convenience,
let $e_{n}^{(k)}=|E(\Gamma_{n}^{(k)})|$.
Then for $n\geq k\geq2$,
the following iteration formula of $e_{n}^{(k)}$ can be obtained by Proposition 2.2 and Theorem 2.3:
$$e_{n}^{(k)}=\mathop{\sum}\limits_{i=1}^{k}e_{n-i}^{(k)}
+\mathop{\sum}\limits_{i=2}^{k}(i-1)F_{n-i+k}^{(k)}.
~~~~~~~~~~~~~~~~~~~~~~~~~~~~~~~~~~~~~~~~~~~~~~~~~~~~~~~~~~~~~~~~~~~~~~~~~~~~~~~~~~~~~~~~~~~~~~~~~~~~~(2.2)$$

By Eqs. (2.1) and (2.2),
the size of $\Gamma_{n}^{(k)}$ can be obtained iteratively.
For example,
the values of $e_{n}^{(k)}$ are listed in Table 2 when $k=2,3,4,5$ and $n=0,1,2,\ldots,14$.

\begin{center}
\includegraphics[scale=0.7715]{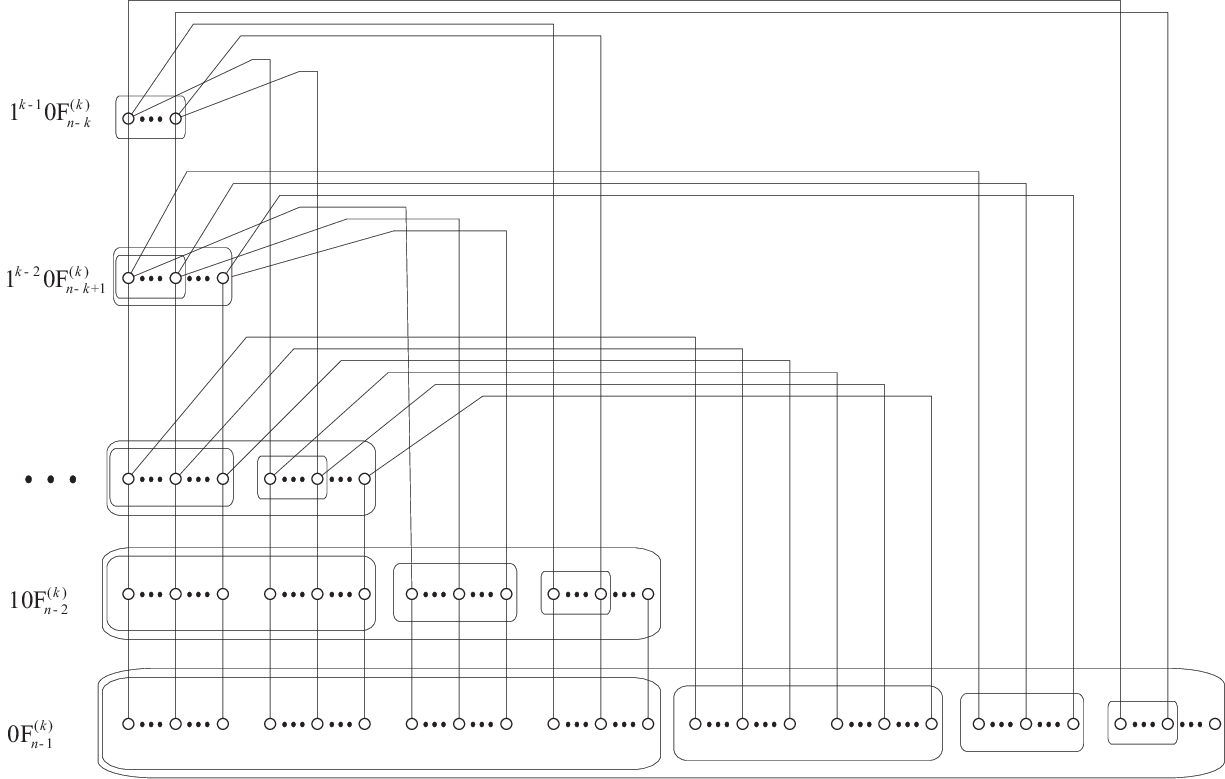}\\
{\footnotesize Fig. 2. Decomposition of the Fibonacci cube $\Gamma_{n}^{(k)}$.}
\end{center}

\begin{table}[!htbp]
\footnotesize
\caption{$e^{(k)}_{n}$, $k=2,3,4,5$ and $n=0,1,2,\ldots,14$.}
\centering
  \begin{tabular}{c|c|c|c|c|c|c|c|c|c|c|c|c|c|c|c}
  \hline
  $n$            & 0& 1 & 2 & 3 & 4 & 5 & 6 & 7 & 8 & 9  & 10 & 11 & 12 & 13 & 14\\
  \hline
   $e^{(2)}_{n}$ & 0& 1 & 2 & 5 & 10& 20&38 & 71&130&235 &420 &744 & 1308 & 2285 & 3970\\
  \hline
  $e^{(3)}_{n}$  & 0&1  &4  &9  &22 &50 &108&230&480&987 &2008&4047& 8094 & 16084 & 31784\\
  \hline
  $e^{(4)}_{n}$  &0&1  &4  &12 &28 &67 &154&344&752&1624&3466&7327& 15368& 32024 & 66356\\
  \hline
  $e^{(5)}_{n}$  & 0&1  &4  &12 &32 &75 &176&402&900&1984&4320&9322& 19956 & 42435 & 89720\\
  \hline
  \end{tabular}
\end{table}

\section{Size of $\Gamma^{(k)}_{n}$}

In this section,
we show that the size of $\Gamma^{(k)}_{n}$ not only can be formulated as the convolution product of
$k$-th order Fibonacci numbers (Subsection 3.1),
but also the linear combination of $k$ consecutive $k$-th order Fibonacci numbers (Subsection 3.2).
Note that in the remainder of this paper,
$\Gamma_{n}^{(2)}, F_{n}^{(2)}$ and $e_{n}^{(2)}$ refer to
$\Gamma_{n}, F_{n}$ and $|E(\Gamma_{n})|$,
respectively.

\subsection{$|E(\Gamma^{(k)}_{n})|$ of the convolution form}

As mentioned before,
$e^{(2)}_{n}$ and $e^{(3)}_{n}$ have been studied in \cite{Klavzar1} and \cite{BelbachirOuldMohamed},
respectively.
Note that in \cite{BelbachirOuldMohamed},
the $k$-th order Fibonacci numbers were introduced in the same way as given by Eq. (2.1),
but the initial values were given by $F^{(3)}_{0}=0$ and $F^{(3)}_{1}=F^{(3)}_{1}=1$.
So there some tiny difference between the formulas of $e^{(3)}_{n}$ given by Theorem 1 of \cite{BelbachirOuldMohamed}
and the result (Theorem 3.2) obtained in this paper.
First we give the following lemma which will be used in the proof of the main result.

\trou \noi {\bf Lemma 3.1}.
\emph{Let $n\geq1$.
Then}

$(i)$ $\mathop{\sum}\limits_{i=1}^{n}i2^{i-1}=(n-1)2^n+1$,

$(ii)$ $\mathop{\sum}\limits_{i=1}^{n}i2^{n-i} + \mathop{\sum}\limits_{i=1}^{n-2}i(n-1-i)2^{n-2-i} =n2^{n-1}$.

\trou \noi {\bf Proof.}
Firstly,
we show that $(i)$ holds by using induction on $n$.
It is obvious that the desired equality holds for $n=1$.
Let $n\geq1$ and assume that it holds for $n$.
Then

$\mathop{\sum}\limits_{i=1}^{n+1}i2^{i-1}$
$=\mathop{\sum}\limits_{i=1}^{n}i2^{i-1}+(n+1)2^{n}$
$=(n-1)2^n+1+(n+1)2^{n}$
$=n2^{n+1}+1$.
\\
This completes the proof of $(i)$.

Now
we show that $(ii)$ holds for $n\geq1$.
For $n=1,2$ the latter sum vanishes,
thus the equality holds for $n=1,2$ by directly computing.
Assume that $n\geq3$ and the equality holds for all indices not more than $n$.
Then by the induction assumption and equality $(i)$,
we have

$\left(\mathop{\sum}\limits_{i=1}^{n+1}i2^{n+1-i} + \mathop{\sum}\limits_{i=1}^{n-1}i(n-i)2^{n-1-i}\right)-n2^{n-1}$

$=\left(\mathop{\sum}\limits_{i=1}^{n+1}i2^{n+1-i} + \mathop{\sum}\limits_{i=1}^{n-1}i(n-i)2^{n-1-i}\right)$
$-\left(\mathop{\sum}\limits_{i=1}^{n}i2^{n-i} + \mathop{\sum}\limits_{i=1}^{n-2}i(n-1-i)2^{n-2-i}\right)$

$=\left(\mathop{\sum}\limits_{i=1}^{n+1}i2^{n+1-i}-\mathop{\sum}\limits_{i=1}^{n}i2^{n-i}\right)$
$+\left(\mathop{\sum}\limits_{i=1}^{n-1}i(n-i)2^{n-1-i}-\mathop{\sum}\limits_{i=1}^{n-2}i(n-1-i)2^{n-2-i}\right)$

$=\mathop{\sum}\limits_{i=1}^{n+1}2^{n+1-i}+\mathop{\sum}\limits_{i=1}^{n-1}(n-i)2^{n-1-i}$
$=\mathop{\sum}\limits_{i=1}^{n+1}2^{i-1}+\mathop{\sum}\limits_{i=1}^{n-1}i2^{i-1}$

$=(2^{n+1}-1)+((n-2)2^{n-1}+1)$

$=n2^{n-1}+2^{n}$.
\\Thus,
we conclude that
$\mathop{\sum}\limits_{i=1}^{n+1}i2^{n+1-i} + \mathop{\sum}\limits_{i=1}^{n-1}i(n-i)2^{n-1-i}$
$=n2^{n-1}+n2^{n-1}+2^{n}=(n+1)2^{n}.$
So we know that $(ii)$ holds for $n+1$.
The proof is completed.
$\Box$

Now we can give the following convolution formula on $e^{(k)}_{n}$.

\trou \noi {\bf Theorem 3.2}.
\emph{Let $k\geq 2$ and $n\geq1$.
Then}
$$e^{(k)}_{n}
=\mathop{\sum}\limits_{j=1}^{k-1}\left(j
\mathop{\sum}\limits_{i=k-1}^{n+k-1-j}F_{i}^{(k)}F_{n+2k-2-j-i}^{(k)}\right).$$

\trou \noi {\bf Proof.}
For a given integer $k\geq2$,
we will prove the equality by induction on $n$:
firstly,
we show the desired equality holds for $1\leq n\leq k$ (Basis),
then assume that the equality is true for all indices smaller than $n$ (Hypothesis),
and finally we show it holds for $n$ (Induction).

Now we begin some preparation for the proof of base case.
For convenience,
let

$S_{n}^{k}
=\mathop{\sum}\limits_{j=1}^{k-1}j\left(
\mathop{\sum}\limits_{i=k-1}^{n+k-1-j}F_{i}^{(k)}F_{n+2k-2-j-i}^{(k)}\right)$
~~~~~~~~~~~~~~~~~~~~~~~~~~~~~~~~~~~~~~~~~~~~~~~~~~~~~~~~~~~~~~~~~~~~~~~~~~~~~$(3.1)$
\\Then we only need to show
$e^{(k)}_{n}=S_{n}^{k}$.


It is know that $F_{k-1}^{(k)}=1$ and
$F_{n}^{(k)}=2^{n-k}$ for $k \leq n \leq 2k-1$ by Eq. (2.1).
As $\Gamma_{n}^{(k)}\cong Q_{n}$ if $1\leq n<k$,
$e_{n}^{(k)}=n2^{n-1}$.
It can be seen that $\Gamma_{k}^{(k)}$ can be obtained by deleting one vertex of $Q_{k}^{(k)}$,
and so $e_{k}^{(k)}=k(2^{k-1}-1)$.
For $n\in[1,k]$,
as the subscript is not more than the superscript in the summation,
Eq. $(3.1)$ can be written as

$S_{n}^{k}
=\mathop{\sum}\limits_{j=1}^{n}\left(j
\mathop{\sum}\limits_{i=k-1}^{n+k-1-j}F_{i}^{(k)}F_{n+2k-2-j-i}^{(k)}\right)$.
~~~~~~~~~~~~~~~~~~~~~~~~~~~~~~~~~~~~~~~~~~~~~~~~~~~~~~~~~~~~~~~$(3.2)$

With these preparations,
we can show that the desired result holds for $n\in[1,k]$ (Basis) as follows.
\\
For $k\geq2$,
$e^{(k)}_{1}=1$ and by Eq. $(3.2)$,
$S_{1}^{k}=F^{(k)}_{k-1}F^{(k)}_{k-1}=1$.
So $e^{(k)}_{n}=S_{n}^{k}$.
\\
For $k=2$,
$e^{(2)}_{2}=2$ and by Eq. $(3.1)$,
$S_{2}^{2}=F^{(2)}_{1}F^{(2)}_{2}+F^{(2)}_{2}F^{(2)}_{1}=2$.
So $e^{(k)}_{n}=S_{n}^{k}$.
\\
For $k\geq3$,
$e^{(k)}_{2}=4$ and by Eq. $(3.2)$,
$S_{2}^{k}=(F^{(k)}_{k-1}F^{(k)}_{k}+F^{(k)}_{k}F^{(k)}_{k-1})+2F^{(k)}_{k-1}F^{(k)}_{k-1}=4$.
So $e^{(k)}_{n}=S_{n}^{k}$.
\\
The remaining case is that $k\geq3$ and $n\geq3$.
It can be seen that this case can be divided into the subcases $k\geq4$ and $3\leq n<k$,
and $k\geq3$ and $n=k$.
For the former,
by Eq. $(3.2)$,
the above preparation and Lemma $3.1(ii)$ we have

$S_{n}^{k}
=\mathop{\sum}\limits_{j=1}^{n}\left(j
\mathop{\sum}\limits_{i=k-1}^{n+k-1-j}F_{i}^{(k)}F_{n+2k-2-j-i}^{(k)}\right)$

$=\mathop{\sum}\limits_{j=1}^{n}\left(jF_{k-1}^{(k)}F_{n+k-1-j}^{(k)}+F_{n+k-1-j}^{(k)}F_{k-1}^{(k)}\right)
+\mathop{\sum}\limits_{j=1}^{n-2}\left(j
\mathop{\sum}\limits_{i=k}^{n+k-2-j}F_{i}^{(k)}F_{n+2k-2-j-i}^{(k)}\right)$

$=\mathop{\sum}\limits_{j=1}^{n}j\left(2^{n-1-j}+2^{n-1-j}\right)
+\mathop{\sum}\limits_{j=1}^{n-2}\left(j\mathop{\sum}\limits_{i=k}^{n+k-2-j}2^{i-k}2^{n+k-2-j-i}\right)$

$=\mathop{\sum}\limits_{j=1}^{n}j2^{n-j}
+\mathop{\sum}\limits_{j=1}^{n-2}\left(j\mathop{\sum}\limits_{i=k}^{n+k-2-j}2^{n-2-j}\right)$

$=\mathop{\sum}\limits_{j=1}^{n}j2^{n-j}
+\mathop{\sum}\limits_{j=1}^{n-2}j\left(n-1-j\right)2^{n-2-j}$

$=n2^{n-1}$

$=e^{(k)}_{n}$.
\\
For the latter case $k\geq3$ and $n=k$,
by Eq. $(3.1)$,
the preparation work and Lemma $3.1(ii)$ we have

$S_{n}^{k}=\mathop{\sum}\limits_{j=1}^{k-1}\left(j
\mathop{\sum}\limits_{i=k-1}^{2k-1-j}F_{i}^{(k)}F_{3k-2-j-i}^{(k)}\right)$

$=\mathop{\sum}\limits_{j=1}^{k-1}\left(j\left(F_{k-1}^{(k)}F_{2k-1-j}^{(k)}+F_{2k-1-j}^{(k)}F_{k-1}^{(k)}\right)\right)
+\mathop{\sum}\limits_{j=1}^{k-2}\left(j
\mathop{\sum}\limits_{i=k}^{2k-2-j}F_{i}^{(k)}F_{3k-2-j-i}^{(k)}\right)$

$=\mathop{\sum}\limits_{j=1}^{k-1}j\left(2^{k-1-j}+2^{k-1-j}\right)
+\mathop{\sum}\limits_{j=1}^{k-2}\left(j\mathop{\sum}\limits_{i=k}^{2k-2-j}2^{i-k}2^{2k-2-j-i}\right)$

$=\mathop{\sum}\limits_{j=1}^{k-1}j2^{k-j}
+\mathop{\sum}\limits_{j=1}^{k-2}\left(j\mathop{\sum}\limits_{i=k}^{2k-2-j}2^{k-2-j}\right)$

$=\mathop{\sum}\limits_{j=1}^{k-1}j2^{k-j}
+\mathop{\sum}\limits_{j=1}^{k-2}j\left(k-j-1\right)2^{k-2-j}$

$=\left(\mathop{\sum}\limits_{j=1}^{k}j2^{k-j}
+\mathop{\sum}\limits_{j=1}^{k-2}j(k-j-1)2^{k-2-j}\right)-k2^{k-k}$

$=k2^{k-1}-k$

$=e_{k}^{k}$.

Let $n(\geq k+1)$ be an integer,
and assume that the desired result holds for all indices smaller than $n$ (Hypothesis).
So for $1\leq s \leq k$,

$e_{n-s}^{(k)}=S_{n-s}^{(k)}=\mathop{\sum}\limits_{j=1}^{k-1}
\left(j\mathop{\sum}\limits_{i=k-1}^{n+k-s-1-j}F_{i}^{(k)}F_{n-s+2k-2-j-i}^{(k)}\right)$.
\\
Note that in the above sum if $i>n-s+k-1-j$,
then

$(n-s+2k-2-j)-i<(n-s+2k-2-j)-(n-s-1+k-j)=k-1$.
\\This means that $F_{n-s+2k-2-j-i}^{(k)}=0$ and so

$e_{n-s}^{(k)}=S_{n-s}^{(k)}=\mathop{\sum}\limits_{j=1}^{k-1}
\left(j\mathop{\sum}\limits_{i=k-1}^{n+k-2-j}F_{i}^{(k)}F_{n-s+2k-2-j-i}^{(k)}\right)$.
~~~~~~~~~~~~~~~~~~~~~~~~~~~~~~~~~~~~~~~~$(3.3)$

Then,
using Eq. (2.2),
the induction assumption,
and Eqs. (2.1) and (3.3),
we compute as follows (Induction):

$e_{n}^{k}=\mathop{\sum}\limits_{s=1}^{k}e_{n-s}^{(k)}
+\mathop{\sum}\limits_{t=2}^{k}(t-1)F_{n-t+k}^{(k)}$

$=\mathop{\sum}\limits_{s=1}^{k}S_{n-s}^{(k)}
+\mathop{\sum}\limits_{t=2}^{k}(t-1)F_{n-t+k}^{(k)}$

$=\mathop{\sum}\limits_{s=1}^{k}
\left(\mathop{\sum}\limits_{j=1}^{k-1}
\left(j\mathop{\sum}\limits_{i=k-1}^{n+k-2-j}F_{i}^{(k)}F_{n-s+2k-2-j-i}^{(k)}\right)\right)
+\mathop{\sum}\limits_{j=1}^{k-1}jF_{n-j-1+k}^{(k)}$

$=\mathop{\sum}\limits_{j=1}^{k-1}
\left(j\mathop{\sum}\limits_{i=k-1}^{n+k-2-j}\left(F_{i}^{(k)}\mathop{\sum}\limits_{s=1}^{k}
F_{n-s+2k-2-j-i}^{(k)}\right)\right)
+\mathop{\sum}\limits_{j=1}^{k-1}jF_{n-j-1+k}^{(k)}$

$=\mathop{\sum}\limits_{j=1}^{k-1}
\left(j\mathop{\sum}\limits_{i=k-1}^{n-2+k-j}F_{i}^{(k)}F_{n+2k-2-j-i}^{(k)}\right)
+\mathop{\sum}\limits_{j=1}^{k-1}jF_{n-j-1+k}^{(k)}F_{k-1}^{(k)}$

$=\mathop{\sum}\limits_{j=1}^{k-1}
\left(j\mathop{\sum}\limits_{i=k-1}^{n-1+k-j}F_{i}^{(k)}F_{n+2k-2-j-i}^{(k)}\right)$

$=S_{n}^{(k)}$.
\\This completes the proof.
$\Box$

By Theorem 3.2,
it can be seen that for $k=2,3,4,5$,
and $n\geq1$:

$e^{(2)}_{n}=\mathop{\sum}\limits_{i=1}^{n}F_{i}^{(2)}F_{n+1-i}^{(2)},$

$e^{(3)}_{n}=\mathop{\sum}\limits_{i=2}^{n+1}F_{i}^{(3)}F_{n+3-i}^{(3)}
+2\mathop{\sum}\limits_{i=2}^{n}F_{i}^{(3)}F_{n+2-i}^{(3)},$

$e^{(4)}_{n}=\mathop{\sum}\limits_{i=3}^{n+2}F_{i}^{(4)}F_{n+5-i}^{(4)}
+2\mathop{\sum}\limits_{i=3}^{n+1}F_{i}^{(4)}F_{n+4-i}^{(4)}
+3\mathop{\sum}\limits_{i=3}^{n}F_{i}^{(4)}F_{n+3-i}^{(4)},$ and

$e^{(5)}_{n}=\mathop{\sum}\limits_{i=4}^{n+3}F_{i}^{(5)}F_{n+7-i}^{(5)}
+2\mathop{\sum}\limits_{i=4}^{n+2}F_{i}^{(5)}F_{n+6-i}^{(5)}
+3\mathop{\sum}\limits_{i=4}^{n+1}F_{i}^{(5)}F_{n+5-i}^{(5)}
+4\mathop{\sum}\limits_{i=4}^{n}F_{i}^{(5)}F_{n+4-i}^{(5)}.$

\subsection{$|E(\Gamma^{(k)}_{n})|$ of the linear combination form}

The linear formula of $e^{(k)}_{n}$ is determined as shown by Theorem 3.6 in this subsection.
To obtain it we need to show that Lemmas 3.3,
3.4 and 3.5 hold.

\trou \noi {\bf Lemma 3.3}.
\emph{Let $n\geq k$.
Then}

$(i)$ for $k\geq4$ and $1\leq t\leq k-2$,

$F_{n+k-1}^{(k)}
=\mathop{\sum}\limits_{j=1}^{k-(t+2)}jF_{n+k-2-j}^{(k)}
+(k-t)\mathop{\sum}\limits_{s=0}^{t-1}F_{n+s}^{(k)}
+\mathop{\sum}\limits_{i=1}^{k-t}(k-i-t+1)F_{n-i}^{(k)}$,

$(ii)$ for $k\geq2$,

$\mathop{\sum}\limits_{i=1}^{k}iF_{n-i}^{(k)}
=(2k-1)F_{n}^{(k)}-F_{n+k-1}^{(k)}+\mathop{\sum}\limits_{j=1}^{k-3}jF_{n+k-2-j}^{(k)}$.

\trou \noi {\bf Proof.}
First we consider the equality in $(i)$.
For convenience,
let

$L(t)=\mathop{\sum}\limits_{j=1}^{k-(t+2)}jF_{n+k-2-j}^{(k)}
+(k-t)\mathop{\sum}\limits_{s=0}^{t-1}F_{n+s}^{(k)}
+\mathop{\sum}\limits_{i=1}^{k-t}(k-i-t+1)F_{n-i}^{(k)}$.
\\Then we need to show $L(k-2)=F_{n+k-1}^{(k)}$ and for $1\leq t\leq k-3$,
$L(t)=L(t+1)$.
In fact,
keeping in mind that the recurrence given by Eq. (2.1),
we compute as follows
\\
$L(k-2)=2\mathop{\sum}\limits_{s=0}^{k-3}F_{n+s}^{(k)}
+\mathop{\sum}\limits_{i=1}^{2}(3-i)F_{n-i}^{(k)}$
\\
$=\mathop{\sum}\limits_{s=-2}^{k-3}F_{n+s}^{(k)}+\mathop{\sum}\limits_{s=-1}^{k-3}F_{n+s}^{(k)}$
\\
$=F_{n+k-2}^{(k)}+\mathop{\sum}\limits_{s=-1}^{k-3}F_{n+s}^{(k)}$
\\
$=F_{n+k-1}^{(k)}$,
\\and
\\
$L(t)=\mathop{\sum}\limits_{j=1}^{k-(t+2)}jF_{n+k-2-j}^{(k)}
+(k-t)\mathop{\sum}\limits_{s=0}^{t-1}F_{n+s}^{(k)}
+\mathop{\sum}\limits_{i=1}^{k-t}(k-i-t+1)F_{n-i}^{(k)}$
\\
$=\mathop{\sum}\limits_{j=1}^{k-(t+2)}jF_{n+k-2-j}^{(k)}
+\left((k-t-1)\mathop{\sum}\limits_{s=0}^{t-1}F_{n+s}^{(k)}+\mathop{\sum}\limits_{s=0}^{t-1}F_{n+s}^{(k)}\right)
+\left(\mathop{\sum}\limits_{i=1}^{k-t}F_{n-i}^{(k)}+\mathop{\sum}\limits_{i=1}^{k-t}(k-i-t)F_{n-i}^{(k)}\right)$
\\$=\mathop{\sum}\limits_{j=1}^{k-(t+2)}jF_{n+k-2-j}^{(k)}
+\left(\mathop{\sum}\limits_{s=0}^{t-1}F_{n+s}^{(k)}+\mathop{\sum}\limits_{i=1}^{k-t}F_{n-i}^{(k)}\right)
+(k-t-1)\mathop{\sum}\limits_{s=0}^{t-1}F_{n+s}^{(k)}+\mathop{\sum}\limits_{i=1}^{k-t}(k-i-t)F_{n-i}^{(k)}$
\\$=\mathop{\sum}\limits_{j=1}^{k-(t+2)}jF_{n+k-2-j}^{(k)}
+F_{n+t}^{(k)}
+(k-t-1)\mathop{\sum}\limits_{s=0}^{t-1}F_{n+s}^{(k)}+\mathop{\sum}\limits_{i=1}^{k-t}(k-i-t)F_{n-i}^{(k)}$
\\$=\mathop{\sum}\limits_{j=1}^{k-(t+3)}jF_{n+k-2-j}^{(k)}
+(k-(t+2))F_{n+t}^{(k)}+F_{n+t}^{(k)}
+(k-t-1)\mathop{\sum}\limits_{s=0}^{t-1}F_{n+s}^{(k)}+\mathop{\sum}\limits_{i=1}^{k-t}(k-i-t)F_{n-i}^{(k)}$
\\$=\mathop{\sum}\limits_{j=1}^{k-(t+3)}jF_{n+k-2-j}^{(k)}
+(k-t-1)F_{n+t}^{(k)}
+(k-t-1)\mathop{\sum}\limits_{s=0}^{t-1}F_{n+s}^{(k)}+\mathop{\sum}\limits_{i=1}^{k-t-1}(k-i-t)F_{n-i}^{(k)}$
\\$=\mathop{\sum}\limits_{j=1}^{k-(t+3)}jF_{n+k-2-j}^{(k)}
+(k-t-1)\mathop{\sum}\limits_{s=0}^{t}F_{n+s}^{(k)}+\mathop{\sum}\limits_{i=1}^{k-t-1}(k-i-t)F_{n-i}^{(k)}$
\\
$=L(t+1)$.

Now we turn to consider $(ii)$.
For $k=2$ and $3$,
the sum in the right-hand side of the desired equality vanishes,
and by direct calculation we have

$F_{n-1}^{(2)}+2F_{n-2}^{(2)}=3F_{n}^{(2)}-F_{n+1}^{(2)}$,
and
$F_{n-1}^{(3)}+2F_{n-2}^{(3)}+3F_{n-3}^{(3)}=5F_{n}^{(3)}-F_{n+2}^{(3)}$.
\\
For $k\geq4$,
by Eq. (2.1) and $(i)$,
we have

$\mathop{\sum}\limits_{i=1}^{k}iF_{n-i}^{(k)}
-\left((2k-1)F_{n}^{(k)}-F_{n+k-1}^{(k)}+\mathop{\sum}\limits_{j=1}^{k-3}jF_{n+k-2-j}^{(k)}\right)$

$=F_{n+k-1}^{(k)}
-\mathop{\sum}\limits_{j=1}^{k-3}jF_{n+k-2-j}^{(k)}-(2k-1)F_{n}^{(k)}+
\left(k\mathop{\sum}\limits_{i=1}^{k}F_{n-i}^{(k)}-\mathop{\sum}\limits_{i=1}^{k}(k-i)F_{n-i}^{(k)}\right)$

$=F_{n+k-1}^{(k)}
-\mathop{\sum}\limits_{j=1}^{k-3}jF_{n+k-2-j}^{(k)}-\left((2k-1)F_{n}^{(k)}-kF_{n}^{(k)}\right)
-\mathop{\sum}\limits_{i=1}^{k-1}(k-i)F_{n-i}^{(k)}$

$=F_{n+k-1}^{(k)}
-\mathop{\sum}\limits_{j=1}^{k-3}jF_{n+k-2-j}^{(k)}-(k-1)F_{n}^{(k)}
-\mathop{\sum}\limits_{i=1}^{k-1}(k-i)F_{n-i}^{(k)}$

$=F_{n+k-1}^{(k)}-L(1)$

$=0$.
\\
This completes the proof.
$\Box$

Recall that $F_{0}^{(k)}=\ldots=F_{k-2}^{(k)}=0$,
$F_{k-1}^{(k)}=1$ and
$F_{n}^{(k)}=2^{n-k}$ for $k \leq n \leq 2k-1$ By Eq. (2.1).
If $1\leq n<k$,
then $\Gamma_{n}^{(k)}\cong Q_{n}$ and so $e_{n}^{(k)}=n2^{n-1}$.
As the graph $\Gamma_{k}^{(k)}$ can be obtained from $Q_{k}$ by deleting the vertex $1^{k}$,
$e_{k}^{(k)}=k2^{k-1}-k$.

\trou \noi {\bf Lemma 3.4}.
\emph{Let $k\geq2$.
Then}

$(k+1)e^{(k)}_{n}-2ke^{(k)}_{n-1}=(n+k)F_{n+k-1}^{(k)}-k(n+1)F_{n-1}^{(k)}$.

\trou \noi {\bf Proof.}
We will prove the above equality by induction on $n$ for a given $k\geq2$.
First we show that it holds for $n\in[k]$ (Basis).
\\For $k\geq2$,
if $n=1$,
then $(k+1)e^{(k)}_{1}-2ke^{(k)}_{0}=k+1$ and $(1+k)F_{k}^{(k)}-k(1+1)F_{k-1}^{(k)}=k+1$;
if $n=k$,
then $(k+1)e^{(k)}_{n}-2ke^{(k)}_{n-1}=(k+1)(k2^{k-1}-k)-2k(k-1)2^{k-2}=k(2^{k}-k-1)$ and
$(n+k)F_{n+k-1}^{(k)}-k(n+1)F_{n-1}^{(k)}=(k+k)2^{k-1}-k(k+1)=k(2^{k}-k-1)$.
\\For $k\geq3$,
if $2\leq n\leq k-1$,
then $(k+1)e^{(k)}_{n}-2ke^{(k)}_{n-1}=(k+1)(n2^{n-1}-k)-2k(n-1)2^{n-2}=(n+k)2^{n-1}$ and
$(n+k)F_{n+k-1}^{(k)}-k(n+1)F_{n-1}^{(k)}=(n+k)2^{n-1}$.

Now let $n (\geq k)$ assume that it holds for all indices not more than $n$ (Hypothesis).
Then,
we have

$\begin{cases}
(k+1)e^{(k)}_{n-1}-2ke^{(k)}_{n-2}=(n+k-1)F_{n+k-2}^{(k)}-knF_{n-2}^{(k)},\\
(k+1)e^{(k)}_{n-2}-2ke^{(k)}_{n-3}=(n+k-2)F_{n+k-3}^{(k)}-k(n-1)F_{n-3}^{(k)},\\
\ldots \\
(k+1)e^{(k)}_{n-(k-1)}-2ke^{(k)}_{n-k}=(n+1)F_{n}^{(k)}-k(n-k+2)F_{n-k}^{(k)}.\\
\end{cases}$
\\
If sum the above $k-1$ equalities together,
then we obtain

$(k+1)\mathop{\sum}\limits_{i=2}^{k}e_{n+1-i}^{(k)}$
$=2k\mathop{\sum}\limits_{i=2}^{k}e_{n-i}^{(k)}
+\mathop{\sum}\limits_{i=2}^{k}(n+i-1)F_{n+i-2}^{(k)}
-k\mathop{\sum}\limits_{i=2}^{k}(n-(i-2))F_{n-i}^{(k)}$.~~~~~~~~~~~~~~~~~~~~~$(3.4)$
\\
By Eq. (2.2),

$\mathop{\sum}\limits_{i=2}^{k}e_{n-i}^{(k)}
=e_{n}^{(k)}-e_{n-1}^{(k)}-\mathop{\sum}\limits_{i=2}^{k}(i-1)F_{n-i+k}^{(k)}$.
~~~~~~~~~~~~~~~~~~~~~~~~~~~~~~~~~~~~~~~~~~~~~~~~~~~~~~~$(3.5)$

Then,
using Eq (2.2),
Eqs. $(3.4)$ and $(3.5)$,
the induction assumption,
Eq (2.1),
and Lemma 3.3,
we can show that the desired result holds for $n+1$ as follows (Induction):
\\
$(k+1)e^{(k)}_{n+1}-2ke^{(k)}_{n}$
\\
$=(k+1)\left(\mathop{\sum}\limits_{i=1}^{k}e_{n+1-i}^{(k)}
+\mathop{\sum}\limits_{i=2}^{k}(i-1)F_{n+1-i+k}^{(k)}\right)-2ke^{(k)}_{n}$
\\
$=-(k-1)e^{(k)}_{n}+(k+1)\mathop{\sum}\limits_{i=2}^{k}e_{n+1-i}^{(k)}
+(k+1)\mathop{\sum}\limits_{i=2}^{k}(i-1)F_{n+1-i+k}^{(k)}$
\\
{\scriptsize
$=-(k-1)e^{(k)}_{n}+\left(2k\mathop{\sum}\limits_{i=2}^{k}e_{n-i}^{(k)}
+\mathop{\sum}\limits_{i=2}^{k}(n+i-1)F_{n+i-2}^{(k)}
-k\mathop{\sum}\limits_{i=2}^{k}(n-(i-2))F_{n-i}^{(k)}\right)$
$+(k+1)\mathop{\sum}\limits_{i=2}^{k}(i-1)F_{n+1-i+k}^{(k)}$}
\\
{\scriptsize$=
-(k-1)e^{(k)}_{n}+2k\left(e_{n}^{(k)}-e_{n-1}^{(k)}-\mathop{\sum}\limits_{i=2}^{k}(i-1)F_{n-i+k}^{(k)}\right)
+\mathop{\sum}\limits_{i=2}^{k}(n+i-1)F_{n+i-2}^{(k)}$
$-k\mathop{\sum}\limits_{i=2}^{k}(n-(i-2))F_{n-i}^{(k)}$}
\\{\scriptsize$~~~+(k+1)\mathop{\sum}\limits_{i=2}^{k}(i-1)F_{n+1-i+k}^{(k)}$}
\\
{\scriptsize$=
\left((k+1)e^{(k)}_{n}-2ke_{n-1}^{(k)}\right)-2k\mathop{\sum}\limits_{i=2}^{k}(i-1)F_{n-i+k}^{(k)}
+\mathop{\sum}\limits_{i=2}^{k}(n+i-1)F_{n+i-2}^{(k)}$
$-k\mathop{\sum}\limits_{i=2}^{k}(n-(i-2))F_{n-i}^{(k)}+(k+1)\mathop{\sum}\limits_{i=2}^{k}(i-1)F_{n+1-i+k}^{(k)}$}
\\
{\scriptsize
$=\left((n+k)F_{n+k-1}^{(k)}-k(n+1)F_{n-1}^{(k)}\right)-2k\mathop{\sum}\limits_{i=2}^{k}(i-1)F_{n-i+k}^{(k)}
+\mathop{\sum}\limits_{i=2}^{k}(n+i-1)F_{n+i-2}^{(k)}$
$-k\mathop{\sum}\limits_{i=2}^{k}(n-(i-2))F_{n-i}^{(k)}$}
\\{\scriptsize$~~~+(k+1)\mathop{\sum}\limits_{i=2}^{k}(i-1)F_{n+1-i+k}^{(k)}$}
\\
{\scriptsize
$=\left((n+k)F_{n+k-1}^{(k)}+\mathop{\sum}\limits_{i=2}^{k}(n+i-1)F_{n+i-2}^{(k)}\right)
-2k\mathop{\sum}\limits_{i=2}^{k}(i-1)F_{n-i+k}^{(k)}$
$-\left(k(n+1)F_{n-1}^{(k)}+k\mathop{\sum}\limits_{i=2}^{k}(n-(i-2))F_{n-i}^{(k)}\right)$}
\\{\scriptsize
$~~~+(k+1)\mathop{\sum}\limits_{i=2}^{k}(i-1)F_{n+1-i+k}^{(k)}$}
\\
{\scriptsize
$=\mathop{\sum}\limits_{i=2}^{k+1}(n+i-1)F_{n+i-2}^{(k)}
+(k+1)\mathop{\sum}\limits_{i=2}^{k}(i-1)F_{n+1-i+k}^{(k)}-2k\mathop{\sum}\limits_{i=2}^{k}(i-1)F_{n-i+k}^{(k)}$
$-k\mathop{\sum}\limits_{i=1}^{k}(n-(i-2))F_{n-i}^{(k)}$}
\\
{\scriptsize
$=\mathop{\sum}\limits_{i=1}^{k}(n+k+1-i)F_{n+k-i}^{(k)}
+\mathop{\sum}\limits_{i=1}^{k-1}(k+1)iF_{n+k-i}^{(k)}-\mathop{\sum}\limits_{i=2}^{k}2k(i-1)F_{n+k-i}^{(k)}$
$-k\mathop{\sum}\limits_{i=1}^{k}(n-(i-2))F_{n-i}^{(k)}$}
\\
{\scriptsize
$=(n+2k+1)F_{n+k-1}^{(k)}+
\mathop{\sum}\limits_{i=2}^{k-1}\left((n+3k+1-ki))F_{n+k-i}^{(k)}\right)$
$+(n+2k+1-2k^2)F_{n}^{(k)}
-k\left(\mathop{\sum}\limits_{i=1}^{k}(n+2)F_{n-i}^{(k)}-\mathop{\sum}\limits_{i=1}^{k}(-i)F_{n-i}^{(k)}\right)$}
\\
{\scriptsize
$=(n+k+1)\mathop{\sum}\limits_{i=1}^{k}F_{n+k-i}^{(k)}
-k\mathop{\sum}\limits_{i=1}^{k}(n+2)F_{n-i}^{(k)}+kF_{n+k-1}^{(k)}
+\mathop{\sum}\limits_{i=2}^{k-1}(2k-ki)F_{n+k-i}^{(k)}$
$-k(2k-1)F_{n}^{(k)}
+k\mathop{\sum}\limits_{i=1}^{k}(-i)F_{n-i}^{(k)}$}
\\
{\scriptsize
$=(n+k+1)F_{n+k}^{(k)}-k(n+2)F_{n}^{(k)}+kF_{n+k-1}^{(k)}
+\mathop{\sum}\limits_{i=2}^{k-1}(2k-ki)F_{n+k-i}^{(k)}-k(2k-1)F_{n}^{(k)}$
$+k\mathop{\sum}\limits_{i=1}^{k}(-i)F_{n-i}^{(k)}$}
\\
{\scriptsize
$=(n+k+1)F_{n+k}^{(k)}-k(n+2)F_{n}^{(k)}$
$+k\left(F_{n+k-1}^{(k)}
-\mathop{\sum}\limits_{j=1}^{k-3}jF_{n+k-2-j}^{(k)}
-(2k-1)F_{n}^{(k)}+\mathop{\sum}\limits_{i=1}^{k}iF_{n-i}^{(k)}\right)$}
\\
$=(n+k+1)F_{n+k-1}^{(k)}-k(n+2)F_{n}^{(k)}$.
\\The proof is completed.
$\Box$

It is know that $F_{n}^{(2)}= 2F_{n+2}^{(2)}-F_{n+3}^{(2)}.$
For general $k\geq2$,
the following parallel result holds.

\trou \noi {\bf Lemma 3.5}.
\emph{Let $k\geq2$ and $n\geq 0$.
Then}

$F_{n}^{(k)}= 2F_{n+k}^{(k)}-F_{n+k+1}^{(k)}.$

\trou \noi {\bf Proof.}
By Eq. (2.1),

$F_{n}^{(k)}=F_{n+k}^{(k)}-F_{n+k-1}^{(k)}-\ldots-F_{n+2}^{(k)}-F_{n+1}^{(k)}$
and
$F_{n+1}^{(k)}=F_{n+k+1}^{(k)}-F_{n+k}^{(k)}-F_{n+k-1}^{(k)}-\ldots-F_{n+2}^{(k)}$.
So

$F_{n}^{(k)}=F_{n+k}^{(k)}-F_{n+k-1}^{(k)}-\ldots-F_{n+2}^{(k)}
-(F_{n+k+1}^{(k)}-F_{n+k}^{(k)}-F_{n+k-1}^{(k)}-\ldots-F_{n+2}^{(k)})$

$=2F_{n+k}^{(k)}-F_{n+k+1}^{(k)}$.
$\Box$

Now we show that for any given integer $k\geq2$,
the calculation formula of $e^{(k)}_{n}$ can be linearly expressed by $k$ consecutive $k$-th order Fibonacci numbers.

\trou \noi {\bf Theorem 3.6}.
\emph{Let $n\geq k\geq2$.
Then}

$e^{(k)}_{n}=\frac{1}{2(2k)^{k}-(k+1)^{k+1}}
\mathop{\sum}\limits_{j=0}^{k-1}(nA_{j}+B_{j})F_{n+j}^{(k)},$
\emph{where}

$A_{j}=(2k)^{k}-k(k+1)^{k}+(k-1)(k+1)^{j}(2k)^{k-j-1},$
\emph{and}

$B_{j}=(2k)^{k-j-1}(k+1)^{j}\left(3k+(k-1)j-1\right)-k(k+1)^{k}$.

\trou \noi {\bf Proof.}
The main idea of the proof of the desired result is as follows:
firstly by Lemma 3.4,
$e^{(k)}_{n}$ is used to replace $e^{(k)}_{n-i}$ ($i\in[1,k]$) in Eq. (2.2);
and then $e^{(k)}_{n}$ is taken as an unknown variable and obtained by solving the equality got above.
The specific process is as follows.
\\
By using this relation shown in Lemma 3.4 repeatedly for $i$ from $1$ to $k$,
we can get the following relation between $e^{(k)}_{n-i}$ and $e^{(k)}_{n}$:

$e^{(k)}_{n-1}=\frac{k+1}{2k}e^{(k)}_{n}+\frac{k(n+1)F_{n-1}^{(k)}-(n+k)F_{n+k-1}^{(k)}}{2k}$,

$e^{(k)}_{n-2}=\frac{k+1}{2k}e^{(k)}_{n-1}+\frac{knF_{n-2}^{(k)}-(n+k-1)F_{n+k-2}^{(k)}}{2k}$

$~~~~~=(\frac{k+1}{2k})^{2}e^{(k)}_{n}
+\frac{k+1}{2k}\frac{k(n+1)F_{n-1}^{(k)}-(n+k)F_{n+k-1}^{(k)}}{2k}
+\frac{knF_{n-2}^{(k)}-(n+k-1)F_{n+k-2}^{(k)}}{2k}$,

$\ldots$

$e^{(k)}_{n-i}=\frac{k+1}{2k}e^{(k)}_{n-(i-1)}+\frac{k(n-i+2)F_{n-i}^{(k)}-(n-i+k+1)F_{n-i+k}^{(k)}}{2k}$

$~~~~~=\ldots$

$~~~~~=(\frac{k+1}{2k})^{i}e^{(k)}_{n}
+\mathop{\sum}\limits_{j=0}^{i-1}
(\frac{k+1}{2k})^{j}\frac{k(n-i+j+2)F_{n-i+j}^{(k)}-(n+k-i+j+1)F_{n+k-i+j}^{(k)}}{2k}$,

$\ldots$

$e^{(k)}_{n-k}=\frac{k+1}{2k}e^{(k)}_{n-(k-1)}+\frac{k(n-k+2)F_{n-k}^{(k)}-(n+1)F_{n}^{(k)}}{2k}$

$~~~~~=\ldots$

$~~~~~=(\frac{k+1}{2k})^{k}e^{(k)}_{n}
+\mathop{\sum}\limits_{j=0}^{k-1}
(\frac{k+1}{2k})^{j}\frac{k(n-k+j+2)F_{n-k+j}^{(k)}-(n+j+1)F_{n+j}^{(k)}}{2k}$.
\\
By Eq. (2.2) and the above relation between $e^{(k)}_{n-i}$ and $e^{(k)}_{n}$,
we have

$e_{n}^{(k)}=\mathop{\sum}\limits_{i=1}^{k}e_{n-i}^{(k)}+\mathop{\sum}\limits_{i=2}^{k}(i-1)F_{n-i+k}^{(k)}$

$~~~~=\mathop{\sum}\limits_{i=1}^{k}\left((\frac{k+1}{2k})^{i}e^{(k)}_{n}
+\mathop{\sum}\limits_{j=0}^{i-1}
(\frac{k+1}{2k})^{j}\frac{k(n-i+j+2)F_{n-i+j}^{(k)}-(n+k-i+j+1)F_{n+k-i+j}^{(k)}}{2k}\right)
+\mathop{\sum}\limits_{i=2}^{k}(i-1)F_{n-i+k}^{(k)}$

$~~~~=\mathop{\sum}\limits_{i=1}^{k}(\frac{k+1}{2k})^{i}e^{(k)}_{n}
+\mathop{\sum}\limits_{i=1}^{k}\mathop{\sum}\limits_{j=0}^{i-1}
(\frac{k+1}{2k})^{j}\frac{k(n-i+j+2)F_{n-i+j}^{(k)}-(n+k-i+j+1)F_{n+k-i+j}^{(k)}}{2k}
+\mathop{\sum}\limits_{i=2}^{k}(i-1)F_{n-i+k}^{(k)}$

$~~~~=\mathop{\sum}\limits_{i=1}^{k}(\frac{k+1}{2k})^{i}e^{(k)}_{n}
+\mathop{\sum}\limits_{s=1}^{k}\mathop{\sum}\limits_{t=0}^{k-s}
(\frac{k+1}{2k})^{t}\frac{k(n-s+2)F_{n-s}^{(k)}-(n+k-s+1)F_{n+k-s}^{(k)}}{2k}
+\mathop{\sum}\limits_{j=1}^{k-1}jF_{n+k-j-1}^{(k)}$.
\\So

$(\mathop{\sum}\limits_{i=1}^{k}(\frac{k+1}{2k})^{i}-1)e^{(k)}_{n}
=\mathop{\sum}\limits_{s=1}^{k}\mathop{\sum}\limits_{t=0}^{k-s}
(\frac{k+1}{2k})^{t}\frac{(n+k-s+1)F_{n+k-s}^{(k)}-k(n-s+2)F_{n-s}^{(k)}}{2k}
-\mathop{\sum}\limits_{j=1}^{k-1}jF_{n+k-j-1}^{(k)}$.
\\Thus,
\\
$(\mathop{\sum}\limits_{i=1}^{k}(k+1)^{i}(2k)^{k-i}-(2k)^{k})e^{(k)}_{n}$
\\
$=\mathop{\sum}\limits_{s=1}^{k}\mathop{\sum}\limits_{t=0}^{k-s}
(k+1)^{t}(2k)^{k-1-t}((n+k-s+1)F_{n+k-s}^{(k)}-k(n-s+2)F_{n-s}^{(k)})
-(2k)^{k}\mathop{\sum}\limits_{j=1}^{k-1}jF_{n+k-j-1}^{(k)}$
\\
$=${\scriptsize $\left((2k)^{k-1}+(k+1)(2k)^{k-2}+\ldots+(k+1)^{k-3}(2k)^{2}+(k+1)^{k-2}(2k)+(k+1)^{k-1}\right)
\left((n+k)F_{n+k-1}^{(k)}-k(n+1)F_{n-1}^{(k)}\right)$}
\\
$~~+${\scriptsize $\left((2k)^{k-1}+(k+1)(2k)^{k-2}+\ldots+(k+1)^{k-3}(2k)^{2}+(k+1)^{k-2}(2k)\right)
\left((n+k-1)F_{n+k-2}^{(k)}-knF_{n-2}^{(k)}\right)-(2k)^{k}F_{n+k-2}^{(k)}$}
\\
$~~+${\scriptsize $\left((2k)^{k-1}+(k+1)(2k)^{k-2}+\ldots+(k+1)^{k-3}(2k)^{2}\right)
\left((n+k-2)F_{n+k-3}^{(k)}-k(n-1)F_{n-3}^{(k)}\right)-(2k)^{k}F_{n+k-3}^{(k)}$}
\\
$\ldots$
\\
$~~+${\scriptsize $\left((2k)^{k-1}+(k+1)(2k)^{k-2}\right)
\left((n+2)F_{n+1}^{(k)}-k(n-k+3)F_{n-k+1}^{(k)}\right)-(k-2)(2k)^{k}F_{n+1}^{(k)}$}
\\
$~~+${\scriptsize $(2k)^{k-1}
\left((n+1)F_{n}^{(k)}-k(n-k+2)F_{n-k}^{(k)}\right)-(k-1)(2k)^{k}F_{n}^{(k)}$}.
~~~~~~~~~~~~~~~~~~~~~~~~~~~~~~~~~~~~~~~~~~~~~~~~~~~~~~~~~$(3.6)$
\\
Note that by Eq. (2.1),

$(n+k)F_{n+k-1}^{(k)}-k(n+1)F_{n-1}^{(k)}$

$=(n+k)F_{n+k-1}^{(k)}-k(n+1)(F_{n+k-1}^{(k)}-F_{n+k-2}^{(k)}-\ldots-F_{n+1}^{(k)}-F_{n}^{(k)})$

$=(1-k)nF_{n+k-1}^{(k)}+k(n+1)(F_{n+k-2}^{(k)}+\ldots+F_{n+1}^{(k)}+F_{n}^{(k)})$;
~~~~~~~~~~~~~~~~~~~~~~~~~~~~~~~~~~~~~~~~~~~~~~~~~~~~~$(3.7)$
\\
and if $2\leq s\leq k$,
then by Lemma 3.5,

$(n+k-s+1)F_{n+k-s}^{(k)}-k(n-s+2)F_{n-s}^{(k)}$

$=(n+k-s+1)F_{n+k-s}^{(k)}-k(n-s+2)(2F_{n+k-s}^{(k)}-F_{n+k-s+1}^{(k)})$

$=\left((2k-1)s+1-3k-(2k-1)n\right)F_{n+k-s}^{(k)}+k(n-s+2)F_{n+k-s+1}^{(k)}$.
~~~~~~~~~~~~~~~~~~~~~~~~~~~~$(3.8)$
\\
Substituting $(3.7)$ and $(3.8)$ into $(3.6)$,
we obtain
\\
$(\mathop{\sum}\limits_{i=1}^{k}(k+1)^{i}(2k)^{k-i}-(2k)^{k})e^{(k)}_{n}$
\\
$=${\scriptsize $n\left((2k)^{k-1}+(k+1)(2k)^{k-2}+\ldots
+(k+1)^{k-3}(2k)^{2}+(k+1)^{k-2}(2k)-(k-1)(k+1)^{k-1}\right)F_{n+k-1}^{(k)}$}
\\
$+${\scriptsize $n\left((2k)^{k-1}+(k+1)(2k)^{k-2}+\ldots
+(k+1)^{k-3}(2k)^{2}-(k-1)(k+1)^{k-2}(2k)+k(k+1)^{k-1}\right)F_{n+k-2}^{(k)}$}
\\
$+${\scriptsize $n\left((2k)^{k-1}+(k+1)(2k)^{k-2}+\ldots
-(k-1)(k+1)^{k-3}(2k)^{2}+k(k+1)^{k-2}(2k)+k(k+1)^{k-1}\right)F_{n+k-3}^{(k)}$}
\\
$\ldots$
\\
$+${\scriptsize $n\left((2k)^{k-1}-(k-1)(k+1)(2k)^{k-2}+\ldots
+k(k+1)^{k-3}(2k)^{2}+k(k+1)^{k-2}(2k)+k(k+1)^{k-1}\right)F_{n+1}^{(k)}$}
\\
$+${\scriptsize $n\left(-(k-1)(2k)^{k-1}+k(2k)^{k-2}+\ldots
+k(k+1)^{k-3}(2k)^{2}+k(k+1)^{k-2}(2k)+k(k+1)^{k-1}\right)F_{n}^{(k)}$}
\\
$+$ {\scriptsize $[0((2k)^{k-1}+(k+1)(2k)^{k-2}+\ldots
+(k+1)^{k-3}(2k)^{2}+(k+1)^{k-2}(2k))+0((k+1)^{k-1}(2k)^{0})-0(2k)^{k}]F_{n+k-1}^{(k)}$}
\\
$+$ {\scriptsize $[(k-1)((2k)^{k-1}+(k+1)(2k)^{k-2}+\ldots
+(k+1)^{k-3}(2k)^{2})+(2k-1)(k+1)^{k-2}(2k)+k(k+1)^{k-1}-(2k)^{k}]F_{n+k-2}^{(k)}$}
\\
$+$ {\scriptsize $[2(k-1)((2k)^{k-1}+(k+1)(2k)^{k-2}+\ldots)
+(2k-1)(k+1)^{k-3}(2k)^{2}+k(k+1)^{k-2}(2k)+k(k+1)^{k-1}-2(2k)^{k}]F_{n+k-3}^{(k)}$}
\\
$\ldots$
\\
$+$ {\scriptsize $[(k-2)(k-1)(2k)^{k-1}+(k-2)(2k-1)(k+1)(2k)^{k-2}+k(+\ldots+(k+1)^{k-2}(2k)+(k+1)^{k-1})-(k-2)(2k)^{k}]F_{n+1}^{(k)}$}
\\
$+$ {\scriptsize
$[(k-1)(2k-1)(2k)^{k-1}+k((k+1)(2k)^{k-2}+\ldots+(k+1)^{k-2}(2k)+(k+1)^{k-1})-(k-1)(2k)^{k}]F_{n}^{(k)}$}.
\\
By the above equality we know that for every $j$ such that $0\leq j\leq k-1$,
the coefficient of $F_{n+j}^{(k)}$ is of the form $nH_{j}+K_{j}$,
where

$H_{j}=\mathop{\sum}\limits_{i=0}^{j}(k+1)^{i}(2k)^{k-i-1}
+k\mathop{\sum}\limits_{i=j}^{k-1}(k+1)^{i}(2k)^{k-i-1}-k(k+1)^{j}(2k)^{k-j-1},$
\emph{and}

$K_{j}=${\scriptsize $(k-j-1)((k-1)$}
$\mathop{\sum}\limits_{i=0}^{j-1}$
{\scriptsize $(k+1)^{i}(2k)^{k-i-1}+(2k-1)(k+1)^{j}(2k)^{k-j-1}-(2k)^{k}$)}
$+k\mathop{\sum}\limits_{i=j+1}^{k-1}$
{\scriptsize $(k+1)^{i}(2k)^{k-i-1}$.}
\\
Now we obtain the following equality by the above discussions:

$(\mathop{\sum}\limits_{i=1}^{k}(k+1)^{i}(2k)^{k-i}-(2k)^{k})e^{(k)}_{n}
=\mathop{\sum}\limits_{j=0}^{k-1}(nH_{j}+K_{j})F_{n+j}^{(k)}$,
~~~~~~~~~~~~~~~~~~~~~~~~~~~~~~~~~~~~~~~~~~$(3.9)$
\\
As

$\mathop{\sum}\limits_{i=1}^{k}(k+1)^{i}(2k)^{k-i}-(2k)^{k}$

$=\frac{(k+1)(2k)^{k-1}(1-(\frac{k+1}{2k})^{k})}{1-\frac{k+1}{2k}}-(2k)^{k}$

$=\frac{2(2k)^{k}-(k+1)^{k+1}}{k-1}$,

$H_{j}=\frac{1}{2k}\left(\frac{(2k)^{k}(1-(\frac{k+1}{2k})^{j+1})}{1-\frac{k+1}{2k}}\right)
+\frac{1}{2k}\left(k\frac{(k+1)^{j+1}(2k-j-1)^{k}(1-(\frac{k+1}{2k})^{k-j-1})}{1-\frac{k+1}{2k}}\right)
-k(k+1)^{j}(2k)^{k-j-1}$

$~~~~=\frac{(2k)^{k}-k(k+1)^{k}+(k-1)(k+1)^{j}(2k)^{k-j-1}}{k-1}$,
and

$K_{j}=$
{\scriptsize
$(k-j-1)(\frac{k-1}{2k}\frac{(2k)^{k}(1-(\frac{k+1}{2k})^{j})}{1-\frac{k+1}{2k}}+(2k-1)(k+1)^{j}(2k)^{k-j-1}-(2k)^{k})
+\frac{1}{2}\frac{(k+1)^{j+1}(2k)^{k-j-1}(1-(\frac{k+1}{2k})^{k-j-1})}{1-\frac{k+1}{2k}}$}

$~~~~=\frac{(3k+(k-1)j-1)(k+1)^{j}(2k)^{k-j-1}-k(k+1)^{k}}{k-1}$.
\\
Let $A_{j}=\frac{k-1}{2(2k)^{k}-(k+1)^{k+1}}H_{j}$ and $B_{j}=\frac{k-1}{2(2k)^{k}-(k+1)^{k+1}}K_{j}$,
$0\leq j\leq k-1$.
Then the desired result is achieved by Eq. $(3.9)$.
$\Box$

For examples,
for $k=2,3,4$ and $5$,

$e^{(2)}_{n}=\frac{nF_{n+1}^{(2)}+2(n+1)F_{n}^{(2)}}{5}$ for $n\geq2$;

$e^{(3)}_{n}=\frac{28nF_{n+2}^{(3)}+12(3n+2)F_{n+1}^{(3)}+48(n+1)F_{n}^{(3)}}{88}$ for $n\geq3$;

$e^{(4)}_{n}
=\frac{219nF_{n+3}^{(4)}+4(61n+25)F_{n+2}^{(4)}+4(41n+55)F_{n+1}^{(4)}+348(n+1)F_{n}^{(4)}}{563}$ for $n\geq4$;
and

$e^{(5)}_{n}
=\frac{518nF_{n+4}^{(5)}+5(109n+27)F_{n+3}^{(5)}+5(118n+63)F_{n+2}^{(5)}
+5(133n+108)F_{n+1}^{(5)}+790(n+1)F_{n}^{(5)}}{1198}$ for $n\geq5$.

\section{Concluding remarks}

The investigation of the properties of Fibonacci cubes attract many researchers \cite{Klavzar,EgeciogluKlavzarMollard},
and also lead to the development of a variety of generalizations and variations.
Among these families of graphs (to list just some of them) are
Fibonacci $(p,r)$-cubes \cite{EgiazarianAstola,Wei1}
Fibonacci-run graphs \cite{EgeciogluIrsic},
Fibonacci $p$-cubes \cite{WeiYang},
Pell graphs \cite{Munarini},
generalized Pell graphs \cite{IrsicKlavzarTan}
weighted Padovan graphs \cite{IrsicChenowetha},
Lucas-run graphs \cite{Wei},
Horadam-Lucas Cubes \cite{TanPodrug} and
associated Mersenne graphs \cite{WeiYang2}.
Based on the decomposition of these graphs,
their sizes were given by iteration formula.
\emph{It would be particularly interesting to find the convolution and linear forms of their sizes.}

Take Fibonacci $p$-cube $\Gamma^{p}_{n}$ \cite{WeiYang,Mollard2} for example,
the iteration and linear formulas of $|E(\Gamma^{p}_{n})|$ have been given.
We close this section by determining the convolution forms of $|E(\Gamma^{p}_{n})|$.

Let $p\geq1$ and $n\geq0$.
Then the \emph{Fibonacci $p$-numbers} $F_{n}^{p}$ are given by the following recurrence relation:
$$F^{p}_{0}=0, F^{p}_{1}=1,\ldots,F^{p}_{p}=1,
F_{n}^{p}=F_{n-1}^{p}+F_{n-p-1}^{p},n>p.
~~~~~~~~~~~~~~~~~~~~~~~~~~~~~~~~~~~~~~~~~~~~~~~~~~~~~~~~~~(4.1)$$

A \emph{Fibonacci $p$-string of length $n$} is a binary string where consecutive 1s are separated by at least $p$ 0s,
and \emph{Fibonacci $p$-cube} $\Gamma^{p}_{n}$ is the subgraph of $Q_{n}$
induced by all the Fibonacci $p$-strings of length $n$.
The case $k=2$ corresponds to the Fibonacci cube $\Gamma_{n}$.
It was shown that the order of $\Gamma^{p}_{n}$ is $F_{n+p+1}^{p}$,
and the size is as follows.

\trou \noi {\bf Theorem 4.1}\cite{WeiYang}.
\emph{Let $p\geq1$ and $n\geq p+1$. Then}

$(i)$ \emph{$|E(\Gamma^{p}_{n})|=|E(\Gamma^{p}_{n-1})|
+|E(\Gamma^{p}_{n-p-1})|+F^{p}_{n}$, and}

$(ii)$ \emph{$|E(\Gamma^{p}_{n})|=\frac{p^{p}nF^{p}_n+
\mathop{\Sigma}\limits^{p}_{t=0}p^{t}(p+1)^{p-t}(n+p-t)F^{p}_{n-t}}{p^{p}+(p+1)^{p+1}}$.}

Now we can show there is a convolution formula for $|E(\Gamma^{p}_{n})|$.

\trou \noi {\bf Theorem 4.2}.
\emph{Let $p\geq1$ and $n\geq p+1$. Then}

\emph{$|E(\Gamma^{p}_{n})|=\mathop{\sum}\limits_{i=1}^{n}F^{p}_{i}F^{p}_{n-i+1}$.}

\trou \noi {\bf Proof.}
Let $p\geq1$ be a given integer.
It can be seen that if $1\leq n \leq p+1$,
then $F^{p}_{n}=1$ by Eq. (4.1),
and $\Gamma^{p}_{n}$ is a star graph,
and so $|E(\Gamma^{p}_{n})|=n$.
For convenience,
let $T^{p}_{n}=\mathop{\sum}\limits_{i=1}^{n}F^{p}_{i}F^{p}_{n-i+1}$.
We will use induction on $n$ to prove that $|E(\Gamma^{p}_{n})|=T^{p}_{n}$.
In fact,
for $1\leq n\leq p+1$,
$T^{p}_{n}=\mathop{\sum}\limits_{i=1}^{n}F^{p}_{i}F^{p}_{n-i+1}=n$,
and so $|E(\Gamma^{p}_{n})|=T^{p}_{n}$ (Basis).
Now we suppose that $|E(\Gamma^{p}_{n})|=T^{p}_{n}$ for all indices smaller than $n$ (Hypothesis).
By using Theorem 4.1 $(i)$,
Eq. (4.1) and induction assumption,
we have (Induction):

$|E(\Gamma^{p}_{n})|$

$=T^{p}_{n-1}+T^{p}_{n-p-1}+F^{p}_{n}$

$=\mathop{\sum}\limits_{i=1}^{n-1}F^{p}_{i}F^{p}_{n-i}
+\mathop{\sum}\limits_{i=1}^{n-p-1}F^{p}_{i}F^{p}_{n-p-i}+F^{p}_{n}$

$=\mathop{\sum}\limits_{i=1}^{n-p-1}F^{p}_{i}(F^{p}_{n-i}+F^{p}_{n-p-i})
+\mathop{\sum}\limits_{i=n-p}^{n-1}F^{p}_{i}F^{p}_{n-i}
+F^{p}_{n}$

$=\mathop{\sum}\limits_{i=1}^{n-p-1}F^{p}_{i}F^{p}_{n-i+1}
+\mathop{\sum}\limits_{i=n-p}^{n-1}F^{p}_{i}F^{p}_{n-i+1}
+F^{p}_{n}F^{p}_{1}$

$=\mathop{\sum}\limits_{i=1}^{n}F^{p}_{i}F^{p}_{n-i+1}$

$=T^{p}_{n}$.
\\
This completes the proof.
$\Box$

\section*{Declaration of competing interest}

The authors declare that they have no known competing financial interests or personal
relationships that could have appeared to influence the work reported in this paper.

\section*{Data availability}

No data was used for the research described in the article.

\end{document}